\documentclass[12pt, reqno]{amsart}
\usepackage{amsmath, amsthm, amscd, amsfonts, amssymb, graphicx, color}
\usepackage[bookmarksnumbered, colorlinks, plainpages]{hyperref}

\newtheorem{theorem}{Theorem}[section]
\newtheorem{lemma}[theorem]{Lemma}

\newtheorem{corollary}[theorem]{Corollary}
\theoremstyle{definition}
\newtheorem{definition}[theorem]{Definition}

\theoremstyle{remark}
\newtheorem{remark}[theorem]{Remark}
\numberwithin{equation}{section}

\begin{document}

\title[Bi-additive $s$-functional inequalities and biderivations in modular spaces]
{Bi-additive $s$-functional inequalities and biderivation in modular spaces}

\author[T. L. Shateri]{T. L. Shateri}
\address{Tayebe Laal Shateri \\ Department of Mathematics and Computer
Sciences, Hakim Sabzevari University, Sabzevar, P.O. Box 397, IRAN}
\email{ \rm t.shateri@hsu.ac.ir; t.shateri@gmail.com}

\thanks{*The corresponding author:
t.shateri@hsu.ac.ir ; t.shateri@gmail.com (Tayebe Laal Shateri)}
 \maketitle

\begin{abstract}
In this article, by using the fixed point method, we prove the generalized
Hyers--Ulam stability of biderivations from an algebra to a modular space, associated to bi-additive $s$-functional inequalities. \vskip 3mm
\end{abstract}

\section{Introduction and preliminaries}\vskip 2mm
Let $\mathcal A$ be an algebra over the real or complex field
$\mathbb{F}$ and let $\mathcal X$ be an $\mathcal A$-module. An
additive mapping $d:\mathcal A \to \mathcal X$ is said to be a
derivation if the functional equation $d(xy)=xd(y)+d(x)y$ holds for
all $x,y\in \mathcal A$. Furthermore, if $d(\alpha x)=\alpha d(x)$
is valid for all $x\in \mathcal A$ and for all $\alpha \in \mathbb
F$, then $d$ is called a linear derivation. A bi-additive mapping $d:\mathcal A\times \mathcal{A}\to \mathcal X$ is called a biderivation, if $d$ is a derivation on each component.

The problem of stability of functional equations was formulated by
 Ulam \cite{Ul} in $1940$, concerning the stability of group homomorphisms. In the year 1941, Hyers \cite{Hy} provided a partial solution to Ulam's problem. Hyers's theorem was
 generalized by Aoki \cite{Ao} for additive mappings and by Rassias \cite{RAS1}
 for linear mappings by considering an unbounded Cauchy difference.
 During the last decades several stability problems for
various functional equations have been investigated by several
authors. We refer the reader to the monographs
\cite{BU,Cz,GA,GV,RAS2}. The stability result concerning
derivations on operator algebras was first obtained by
$\check{S}$emrl \cite{SE}. The reader
may found more about the stability of derivations in \cite{EG,
EM,NI,P1,P2}. Shateri and Sadeghi in \cite{SH} investigated the
generalized Hyers--Ulam stability and superstability for derivations from an  algebra to a modular space.
The theory of modular spaces were founded by Nakano \cite{NA} and
were intensively developed by Luxemburg \cite{L}, Koshi and
Shimogaki \cite{KS} and Yamamuro \cite{YA} and their collaborators.
In the present time the theory of modulars and modular spaces is
extensively applied, in particular, in the study of various Orlicz
spaces \cite{O} and interpolation theory \cite{KR, MA,M}.
\begin{definition}\label{D}
Let $\mathcal{X}$ be an arbitrary vector space.\\
$(a)$ A functional $\rho: \mathcal{X}\to [0,\infty]$ is called a modular if for arbitrary $x, y \in \mathcal{X}$,\\
 $(i) \;\rho(x)=0$ if and only if $x=0$,\\
 $(ii) \;\rho(\alpha x)=\rho(x)$ for every scaler $\alpha$ with $|\alpha|=1$,\\
 $(iii) \;\rho(\alpha x+ \beta y)\leq \rho(x)+\rho(y)$ if and only if $\alpha+\beta=1$ and $\alpha,\beta\geq 0$,\\
$(b)$ if $(iii)$ is replaced by \\
 $(iii)^{'}\; \rho(\alpha x+ \beta y)\leq \alpha \rho(x)+\beta \rho(y)$ if and only if $\alpha+\beta=1$ and $\alpha,\beta\geq 0$,\\
then we say that $\rho$ is a convex modular.
\end{definition}
A modular $\rho$ defines a corresponding modular space, i.e., the
vector space $\mathcal{X}_{\rho}$ given by
\begin{equation*}
\mathcal{X}_{\rho}=\left\{x \in \mathcal{X} :\quad \rho(\lambda
x)\to 0  \mbox{ as }  \lambda \to 0\right\}.
\end{equation*}
Let $\rho$ be a convex modular, the modular space
$\mathcal{X}_{\rho}$ can be equipped with a norm called the
Luxemburg norm, defined by
\begin{equation*}
\|x\|_{\rho}=\inf\left\{\lambda>0 \quad ;\quad \rho\left(\frac{x}{\lambda}\right)\leq 1\right\}.
\end{equation*}

A function modular is said to satisfy the $\Delta_2$--condition if
there exists $\kappa>0$ such that $\rho(2x)\leq \kappa\rho(x)$ for
all $x\in \mathcal{X}_{\rho}$.
\begin{definition}
Let $\{x_n\}$ and $x$ be in $\mathcal{X}_{\rho}$. Then\\
(i) the sequence $\{x_n\}$, with $x_n \in \mathcal{X}_{\rho}$, is $\rho$--convergent to $x$
and write $x_n\stackrel{\rho}{\longrightarrow}x$ if $\rho(x_n-x)\to 0$ as $n \to \infty$.\\
(ii) The sequence  $\{x_n\}$, with $x_n \in \mathcal{X}_{\rho}$, is called $\rho$--Cauchy if $\rho(x_n-x_m)\to 0$ as $n,m \to \infty$.\\
(iii) A subset $\mathcal{S}$ of $\mathcal{X}_{\rho}$ is called
$\rho$--complete complete if and only if any $\rho$--Cauchy
sequence is $\rho$--convergent to an element of $ \mathcal{S}$.
\end{definition}

The modular $\rho$ has the Fatou property if and only if
$\rho(x)\leq \liminf_{n\to\infty}\rho(x_n)$ whenever the sequence
$\{x_n\}$ is $\rho$--convergent to $x$.
\begin{remark}
Note that $\rho(.x)$ is an increasing function, for all $x\in\mathcal X$. Suppose $0<a<b$, then property (iii) of Definition \ref{D} with $y=0$ shows that
$\rho(ax)=\rho\left(\frac{a}{b}bx\right)\leq \rho(bx)$ for all $x\in\mathcal{X}$.
Moreover, if $\rho$ is a convex modular on $\mathcal{X}$ and
$|\alpha|\leq1$, then $\rho(\alpha x)\leq \alpha\rho(x)$  and also
$\rho(x)\leq\frac{1}{2}\rho(2x)$ for all $x\in\mathcal{X}$.
\end{remark}
Let $\mathcal A$ be an algebra and let $f:\mathcal A\times \mathcal A\to {\mathcal X}_{\rho}$ be a mapping. In this paper, by using some ideas of \cite{PA2,PA3} 
we solve the bi-additive $s$-functional inequalities 
\begin{align}\label{eq0}
\rho\big(f(\lambda (x+y),z+w)&+f(\lambda (x+y),z-w)+f(\lambda (x-y),z+w)\nonumber\\&+f(\lambda (x-y),z-w)-4\lambda f(x,z)\big)\nonumber\\
&\leq \rho\Big(4s\big[f(\frac{x+y}{2},z-w)+f(\frac{x-y}{2},z+w)-f( x,z)+f(y,w)\big]\Big)
\end{align}
and
\begin{align}\label{eq00}
 &\rho\Big(4\big[f(\lambda\frac{x+y}{2},z-w)+f(\lambda\frac{x-y}{2},z+w)-\lambda f(x,z)-\lambda f(y,w)\big]\Big)\nonumber\\
 &\leq\rho\Big(s\left[f(x+y,z+w)+f(x+y,z-w)+f(x-y,z+w)+f(x-y,z-w)-4f(x,z)\right]\Big)
\end{align}
where s is a fixed nonzero complex number with $|s|<1$. Then we prove the
generalized Hyers--Ulam stability for biderivations from an algebra to a  modular space, associated to above bi-additive $s$-functional inequalities.
\section{\bf Main results}

 Throughout this paper, we assume that $\mathcal A$ be an algebra and $\mathcal X$ a $\rho$--complete modular space with the convex modular $\rho$ in which has the Fatou property and satisfies the $\Delta_2$--condition with $0<\kappa\leq2$. 
 
\begin{lemma}\label{lem1}\textbf{}
Let $f:\mathcal A\times \mathcal A\to\mathcal X_{\rho}$ be a mapping with
$f(x,0)=f(0,x)=0$ and satisfy in the inequality of the form
\begin{align}\label{eq1}
\rho\Big(f(\lambda (x+y),z+w)&+f(\lambda (x+y),z-w)+f(\lambda (x-y),z+w)\nonumber\\&+f(\lambda (x-y),z-w)-4\lambda f(x,z)\Big)\nonumber\\
&\leq \rho\Big(4s\big[f(\frac{x+y}{2},z-w)+f(\frac{x-y}{2},z+w)-f( x,z)+f(y,w)\big]\Big)
\end{align}
for all $\lambda \in \mathbb T,s<1$ and $x,y,z,w \in \mathcal A$. Then $f$ is biadditive.
\begin{proof}
Letting $x=y,w=0$ and $\lambda=1$ in (\ref{eq1}), we obtain $2f(2x,z)=4f(x,z)$, for all $x,z\in\mathcal X$. So we have
\begin{align}\label{eq2}
&\rho\left(f(x+y,z+w)+f(x+y,z-w)+f(x-y,z+w)+f(x-y,z-w)-4f(x,z)\right)\nonumber\\
&\leq \rho \big(4s\big[f(\frac{x+y}{2},z-w)+f(\frac{x-y}{2},z+w)-f(x,z)+f(y,w)\big]\big)\nonumber\\&=\rho\big(s\big[2f(x+y,z-w)+2f(x-y,z+w)-4f(x,z)+4f(y,w)\big]\big)
\end{align}
for all $x,y,z,w\in\mathcal X$. Letting $w=0$ in (\ref{eq2}), since $|s|<1$ we get 
$$f(x+y,z)+f(x-y,z)=2f(x,z),$$ and so $$f(x_1,z)+f(y_1,z)=2f(\frac{x_1+y_1}{2},z)=f(x_1+y_1,z),$$ for all $x_1=x+y, y_1=x-y$ and $z\in \mathcal A$. Since $|s|<1$ and $f(0,z)=0$, for all $z\in \mathcal A$, we deduce that $f$ is additive in the first variable. Similarly, one can show that $f$ is additive in the second variable. Hence
$f$ is bi-additive.
\end{proof}
\end{lemma}
In following, we give the stability of biderivations from an algebra to a modular space associated to the $s$-functional inequaltiy \ref{eq0}.
\begin{theorem}\label{th.1}
Let $d:\mathcal A\times \mathcal A\to\mathcal X_{\rho}$ satisfy in 
$d(x,0)=d(0,x)=0$ and the inequality of the form
\begin{align}\label{eq1.1}
&\rho\Big(d(\lambda (x+y),z+w)+d(\lambda (x+y),z-w)+d(\lambda (x-y),z+w)\nonumber\\&+d(\lambda (x-y),z-w)-4\lambda d(x,z)\Big)\nonumber\\
&\leq \rho\Big(4s\big[d(\frac{x+y}{2},z-w)+d(\frac{x-y}{2},z+w)-d( x,z)+d(y,w)\big]\Big)\nonumber\\&+\psi(x,y)\psi(z,w)
\end{align}
for all $\lambda \in \mathbb T, s<1$ and $x,y,z,w \in \mathcal A$, where $\psi:\mathcal
A^2\to [0, \infty)$
is a given mapping such that
\begin{equation*}
\psi(2x,2x)\leq 2L \psi(x,x)
\end{equation*}
and
\begin{equation}\label{eq1.2}
\lim_{n\to \infty}\frac{\psi(2^nx,2^ny)}{2^n}=0
\end{equation}
for all $x,y \in \mathcal A$ and a constant $0<L<1$. Then there
exists a unique bi-additive mapping $D:\mathcal A\times \mathcal A\to \mathcal X_{\rho}$ which is $\mathbb C$-linear in the first variable, such that
\begin{equation}\label{eq1.3}
\rho(D(x,z)-d(x,z))\leq \frac{1}{2(1-L)}\psi(x,x)\psi(z,0)
\end{equation}
for all $x,z \in \mathcal A$.

In addition, if the mapping $d:\mathcal A\times \mathcal A\to\mathcal X_{\rho}$ satisfies $d(2x,z)=2d(x,z)$ and
\begin{equation}\label{eq1.4}
\rho \big(d(xy,z)-d(x,z)y-xd(y,z)\big)\leq \psi(x,y)\psi(z,w),
\end{equation}
for all $x,y,z \in \mathcal A$, then $d$ is a biderivation.
\begin{proof}
Consider the set
$$\mathcal E=\{\delta:\mathcal A\times \mathcal A\to \mathcal X_\rho, \quad \delta(0,0)=0\}.$$
We define the function $\widetilde{\rho}$ on $\mathcal E$ as
follows,
\begin{equation*}
\widetilde{\rho}(\delta)=\inf\{c>0 : \rho(\delta(x,z))\leq c \psi(x,x)\psi(y,0)\}.
\end{equation*}
Then $\widetilde{\rho}$ is convex modular. It is sufficient to show
that $\widetilde{\rho}$ satisfies the following condition
$$\widetilde{\rho}(\alpha \delta+\beta \gamma)\leq \alpha\widetilde{\rho}(\delta)+\beta\widetilde{\rho}(\gamma)$$
if $\alpha, \beta\geq0$ such that $\alpha+\beta=1$. Given $\varepsilon>0$, then there exist $c_1>0$ and $c_2>0$ such that
$$c_1\leq \widetilde{\rho}(\delta)+\varepsilon,\quad \rho(\delta(x,z))\leq c_1\psi(x,x)\psi(z,0)$$
and
$$c_2\leq \widetilde{\rho}(\gamma)+\varepsilon,\quad \rho(\gamma(x,z))\leq c_2\psi(x,x)\psi(z,0).$$
If $\alpha,\beta\geq 0$ and $\alpha+\beta=1$, we get
\begin{equation*}
\rho\left(\alpha \delta(x,z)+\beta \gamma(x,z)\right)\leq \alpha
\rho(\delta(x,z))+\beta \rho(\gamma(x,z))\leq (\alpha c_1+\beta
c_2)\psi(x,x)\psi(y,0),
\end{equation*}
hence
\begin{equation*}
\widetilde{\rho}(\alpha \delta+\beta \gamma)\leq \alpha
\widetilde{\rho}(\delta)+\beta \widetilde{\rho}(\gamma)+
(\alpha+\beta)\varepsilon.
\end{equation*}
Therefore
\begin{equation*}
\widetilde{\rho}(\alpha \delta+\beta \gamma)\leq \alpha
\widetilde{\rho}(\delta)+\beta \widetilde{\rho}(\gamma).
\end{equation*}
Moreover, $\widetilde{\rho}$ satisfies the $\Delta_2$--condition
with $0<\kappa<2$. Let $\{\delta_n\}$ be a
$\widetilde{\rho}$--Cauchy sequence in $\mathcal
E_{\widetilde{\rho}}$ and let $\varepsilon>0$ be given. There exists
a positive integer $n_0\in \mathbb{N}$ such that
$\widetilde{\rho}(\delta_n-\delta_m)\leq\varepsilon$ for all
$n,m\geq n_0$.
 Now by considering the definition of the modular $\widetilde{\rho}$, we see that
\begin{equation}\label{eq1.5}
\rho\left(\delta_n(x,z)-\delta_m(x,z)\right)\leq \varepsilon \psi(x,x)\psi(z,0),
\end{equation}
for all $x,z\in\mathcal A$ and $n, m \geq n_0$. Let $x$ be a point of
$\mathcal X$, (\ref{eq1.5}) implies that $\{\delta_n(x,z)\}$ is a
$\rho$--Cauchy sequence in $\mathcal X_{\rho}$. Since $\mathcal
X_{\rho}$ is $\rho$--complete, so $\{\delta_n(x,z)\}$ is
$\rho$--convergent in $\mathcal X_{\rho}$, for each $x,z\in \mathcal
A$. Therefore we can define a function $\delta:\mathcal A\times \mathcal A\to
\mathcal X_{\rho}$ by
$$\delta(x,z)=\lim_{n\to \infty}\delta_n(x,z)$$
for any $x,z\in\mathcal A$. Letting $m\to \infty$, then (\ref{eq1.5})
implies that
 $$\widetilde{\rho}(\delta_n-\delta)\leq \varepsilon $$
for all $n \geq n_0$. Since $\rho$ has the Fatou property, thus
$\{\delta_n\}$ is $\widetilde{\rho}$--convergent sequence in
$\mathcal E_{\widetilde{\rho}}$. consequently $\mathcal
E_{\widetilde{\rho}}$ is  $\widetilde{\rho}$-complete.

Now, we define the function $\mathcal{T}:\mathcal
E_{\widetilde{\rho}}\to\mathcal E_{\widetilde{\rho}}$ as follows
\begin{equation*}
\mathcal{T}\delta(x,z):=\frac{1}{2}\delta(2x,2z)
\end{equation*}
for all $\delta \in \mathcal E_{\widetilde{\rho}}$. Let $\delta,
\gamma \in \mathcal E_{\widetilde{\rho}}$ and let $c\in[0,\infty]$
be an arbitrary constant with $\widetilde{\rho}(\delta-\gamma)\leq
c$. We have
\begin{equation*}
\rho(\delta(x,z)-\gamma(x,z))\leq c \psi(x,x)\psi(z,0)
\end{equation*}
for all $x,z\in \mathcal A$. By the assumption and the last
inequality, we get
\begin{equation*}
\begin{split}
\rho\left(\frac{\delta(2x,2z)}{2}-\frac{\gamma(2x,2z)}{2}\right)&\leq
\frac{1}{2} \rho(\delta(2x,2z)-\gamma(2x,2z))\\&\leq \frac{1}{2}c \psi(2x,2x)\psi(2z,0)\leq Lc \psi(x,x)\psi(z,0)
\end{split}
\end{equation*}
for all $x,z\in \mathcal A$. Hence,
$\widetilde{\rho}(\mathcal{T}\delta-\mathcal{T}\gamma)\leq
L\widetilde{\rho}(\delta-\gamma)$, for all $\delta, \gamma \in
\mathcal E_{\widetilde{\rho}}$ so $\mathcal{T}$ is a
$\widetilde{\rho}$--strict contraction. We show that the
$\widetilde{\rho}$--strict mapping $\mathcal{T}$ satisfies the
conditions of \cite[Theorem 3.4]{Kh}.

Letting $x=y,\lambda=1$ and $w=0$ in (\ref{eq1.1}), and since $\rho(.x)$ is increasing we get
\begin{equation}\label{eq1.6}
\rho\left(d(2x,z)-2d(x,z)\right)\leq\rho\left(2d(2x,z)-4d(x,z)\right)\leq \psi(x,x)\psi(z,0)
\end{equation}
for all $x,z\in\mathcal A$. Replacing $x$ by $2x$ in (\ref{eq1.6}) we
have
\begin{equation*}\label{eq1.7}
\rho\left(d(4x,z)-2d(2x,z)\right)\leq \psi(2x,2x)\psi(z,0)
\end{equation*}
for all $x,z\in\mathcal A$. Since $\rho$ is convex modular, for all $x,z\in\mathcal A$ we
have
\begin{equation*}
\begin{split}
\rho\left(\frac{d(4x,z)}{2}-2d(x,z)\right)&\leq  \frac{1}{2}\rho\big(d(4x,z)-2d(2x,z)\big)+\frac{1}{2}\rho\big(2d(2x,z)-4d(x,z)\big)\\
&\leq \frac{1}{2}\psi(2x,2x)\psi(z,0)+\frac{1}{2}\psi(x,x)\psi(z,0).
\end{split}
\end{equation*}
Moreover,
\begin{equation*}
\begin{split}
\rho\left(\frac{d(2^2x,z)}{2^2}-d(x,z)\right)&\leq
\frac{1}{2}\rho\left(2\frac{d(4x,z)}{2^2}-2d(x,z)\right)\\
&\leq
\frac{1}{2^2}\psi(2x,2x)\psi(z,0)+\frac{1}{2^2}\psi(x,x)\psi(z,0).
\end{split}
\end{equation*}
for all $x,z\in\mathcal A$. By induction we obtain
\begin{equation}\label{eq1.8}
\rho\left(\frac{d(2^nx,z)}{2^n}-d(x,z)\right)\leq
\frac{1}{2^n}\sum_{i=1}^n \psi(2^{i-1}x,2^{i-1}x)\psi(z,0)
\leq \frac{1}{2(1-L)}\psi(x,x)\psi(z,0)
\end{equation}
for all $x,z\in\mathcal A$. Now we assert that
$\delta_{\widetilde{\rho}}(d)=\sup\left\{\widetilde{\rho}\left(\mathcal{T}^n(d)-\mathcal{T}^m(d)\right);
n,m\in \mathbb{N}\right)\}<\infty$. It follows from inequality
(\ref{eq1.6}) that
\begin{align*}
\rho\left(\frac{d(2^nx,z)}{2^n}-\frac{d(2^mx,z)}{2^m}\right)&\leq \frac{1}{2}\rho\left(2\frac{d(2^nx,z)}{2^n}-2d(x,z)
\right)+\frac{1}{2}\rho\left(2\frac{d(2^mx,z)}{2^m}-2d(x,z)\right)\\
&\leq\frac{\kappa}{2}\rho\left(\frac{d(2^nx,z)}{2^n}-d(x,z)\right)+\frac{\kappa}{2}\rho
\left(\frac{d(2^mx,z)}{2^m}-d(x,z)\right)\\
&\leq \frac{1}{1-L}\psi(x,x)\psi(z,0),
\end{align*}
for every $x,z\in\mathcal A$ and $n, m\in \mathbb{N}$, which implies
that
\begin{equation*}\label{eq1.9}
\widetilde{\rho}\left(\mathcal{T}^n(d)-\mathcal{T}^m(d)\right)\leq\frac{1}{1-L},
\end{equation*}
for all $n, m\in \mathbb{N}$. Therefore
$\delta_{\widetilde{\rho}}(d)<\infty$. By \cite[Lemma 3.3]{Kh} we deduce that  $\{\mathcal{T}^n(d)\}$ is $\widetilde{\rho}$-converges to
$D\in\mathcal E_{\widetilde{\rho}}$. Since $\rho$ has the Fatou
property (\ref{eq1.8}) gives
$\widetilde{\rho}(\mathcal{T}D-d)<\infty$.

If we replace $x$ by $2^nx$ in inequality (\ref{eq1.6}), then
\begin{equation*}
\widetilde{\rho}\left(d(2^{n+1}x,z)-2d(2^nx,z)\right)\leq \psi(2^nx,
2^nx)\psi(z,0),
\end{equation*}
for all $x,z\in\mathcal A$. Hence
\begin{equation*}
\begin{split}
\rho\left(\frac{d(2^{n+1}x,z)}{2^{n+1}}-\frac{d(2^nx,z)}{2^n}\right)&\leq
\frac{1}{2^{n+1}}\rho\left(d(2^{n+1}x,z)-2d(2^nx,z)\right)
\leq\frac{1}{2^{n+1}}\psi(2^nx,2^nx)\psi(z,0)\\
&\leq\frac{1}{2^{n+1}}2^n L^n\psi(x,x)\psi(z,0)\leq
\frac{L^n}{2}\psi(x,x)\psi(z,0)\leq\psi(x,x)\psi(z,0)
\end{split}
\end{equation*}
for all $x,z\in\mathcal A$. Therefore
$\widetilde{\rho}(\mathcal{T}(D)-D)<\infty$. It follows from
\cite[Theorem 3.4]{Kh} that $\widetilde{\rho}$--limit of
$\{\mathcal{T}^n(d)\}$ is fixed point of $\mathcal{T}$. If we
replace $x$ by $2^nx$ and $y$ by $2^ny$ and $\lambda=1$ in inequality
(\ref{eq1.1}), then we obtain
\begin{align*}
\rho\Big(\frac{d(2^n(x+y),z+w)}{2^n}&+\frac{d(2^n(x+y),z-w)}{2^n}\\&+\frac{d(2^n(x-y),z+w)}{2^n}
+\frac{d(2^n(x-y),z-w)}{2^n}-4\frac{d(2^nx,z)}{2^n}\Big)\\
&\leq \frac{1}{2^n}\rho\big(d(2^n(x+y),z+w)+d(2^n(x+y),z-w)
\\&+d(2^n(x-y),z+w)+d(2^n(xy),z-w)-4d(2^nx,z)\big)\\
&\leq \rho\Big(\frac{4s}{2^n}\Big[d(2^n(\frac{x+y}{2}),z-w)+d(2^n(\frac{x-y}{2}),z+w)\\&-d(2^nx,z)+d(2^ny,w)\Big]\Big
)+\frac{1}{2^n}\psi(2^nx,2^ny)\psi(z,w)
\end{align*}
for all $x,y,z,w\in\mathcal A$. Hence,
\begin{equation}\label{eq11}
\begin{split}
\rho\big(D(x+y,z+w)&+D(x+y,z-w)+D(x-y,z+w)+D(x-y,z-w)-4D(x,z)\big)\\
&\leq \rho\left(4s\left[D(\frac{x+y}{2},z-w)+D(\frac{x-y}{2},z+w)-D(x,z)+D(y,w)\right]\right)
\end{split}
\end{equation}
for all $x,y,z,w\in\mathcal A$. By Lemma \ref{lem1}, $D$ is biadditive. Put $y=x$ and $w=0$ in (\ref{eq1.1}), we get
$$\rho\big(d(2\lambda x,z)-2\lambda d(x,z)\big)\leq\rho\big(2d(2\lambda x,z)-4\lambda d(x,z)\big)\leq \psi(x,x)\psi(z,0),$$
for all $x,z\in\mathcal A$ and $\lambda \in \mathbb T$. Thus
\begin{equation*}
\rho\Big(\frac{d(2\lambda 2^nx,z)}{2^n}-2\lambda d(2^nx,z)\Big)\leq \frac{1}{2^n}\psi(2^nx,2^nx)\psi(z,0),
\end{equation*}
and so
$$D(2\lambda x,z)-2\lambda D(x,z)=0,$$
for all $x,z\in\mathcal A$ and $\lambda \in \mathbb T$. Therefore 
$D(2\lambda x,z)=2\lambda D(x,z)=0$ and hence $D(\lambda x,z)=\lambda D(x,z)=0$ for all $x,z\in\mathcal A$ and $\lambda \in \mathbb T$. 
Now let $0\neq\lambda\in\mathbb C$ and $M$ an integer greater than $4|\lambda|$. Then $|\frac{\lambda}{M}|<\frac{1}{4}<1-\frac{2}{3}=\frac{1}{3}$. 
By \cite[Theorem 1]{KA1}, there exist three elements $\mu_1,\mu_2,\mu_3\in\mathbb T$ such
that $\frac{3\lambda}{M}=\mu_1+\mu_2+\mu_3$, and $D(x,z) = D(3\frac{1}{3}x,z)=3D(\frac{1}{3}x,z)$, for all $x\in \mathcal A$. So $D(\frac{1}{3}x,z)=\frac{1}{3}D(x,z)$, for all $x\in \mathcal A$. Thus 
\begin{align*}
D(\lambda x,x)&=D\left(\frac{M}{3}.3\frac{\lambda}{M}x,z\right)=M D\left(\frac{1}{3}.3\frac{\lambda}{M}x,z\right)
=\frac{M}{3}D\left(3\frac{\lambda}{M}x\right)\\
&=\frac{M}{3}D\left((\mu_1+\mu_2+\mu_3)x,z\right)=\frac{M}{3}(\mu_1+\mu_2+\mu_3)D(x,z)\\&=\frac{M}{3}.3\frac{\lambda}{M}D(x,z)=\lambda D(x,z),
\end{align*}
for all $x,y,z,w\in\mathcal A$. Hence $D$ is $\mathbb C$-linear in the first variable. Moreover by (\ref{eq1.8}) the inequality (\ref{eq1.3}) is valid. Now if $D^*$ is another fixed point of $\mathcal{T}$, then
\begin{equation*}
\begin{split}
\widetilde{\rho}(D-D^*)&\leq \frac{1}{2}\widetilde{\rho}(2\mathcal{T}(D)-2d)+\frac{1}{2}
\widetilde{\rho}(2\mathcal{T}(D^*)-2d)\\
&\leq
\frac{\kappa}{2}\widetilde{\rho}(\mathcal{T}(D)-d)+\frac{\kappa}{2}\widetilde{\rho}(\mathcal{T}(D^*)-d)
\leq\frac{\kappa}{2(1-L)}<\infty.
\end{split}
\end{equation*}
Since $\mathcal{T}$ is $\widetilde{\rho}$--strict contraction, we get
$$\widetilde{\rho}(D-D^*)=\widetilde{\rho}(\mathcal{T}(D)-\mathcal{T}(D^*))\leq L \widetilde{\rho}(D-D^*),$$
which implies that $\widetilde{\rho}(D-D^*)=0$ or  $D=D^*$, since
$\widetilde{\rho}(D-D^*)<\infty$. This proves the uniqueness of $D$. Now if $d(2x,z)=2d(x,z)$, for all $x,z\in\mathcal A$, then $D(x,z)=\rho-~\lim _{n\to \infty}\frac{d(2^nx,z)}{2^n}=d(x,z)$, for all $x,z\in\mathcal A$. Replacing $x,y$ with $2^nx$ and $2^ny$ in (\ref{eq1.4}) we obtain
\begin{equation*}
\rho \Big(\frac{1}{2^{2n}}\big[d(2^{2n}xy,z)-d(2^nx,z)2^ny-2^nxd(2^ny,z)\big]\Big)
\leq \frac{1}{2^{2n}}\psi(2^nx,2^ny)\psi(z,w),
\end{equation*}
hence $$D(xy,z)-D(x,z)y-xD(y,z)=0$$
for all $x,y,z\in\mathcal A$.  Since $d(x,z)=D(x,z)$, for all $x,z\in\mathcal A$, it follows that $d$ is a biderivation.
\end{proof}
\end{theorem}
Let $\mathcal A$ be a normed algebra. It is known that every normed space is modular space with the modular $\rho(x)=\|x\|$ and $\kappa=2$. If in Theorem \ref{th.1} we put $\psi(x,y)=\sqrt{\theta}(\|x\|^p+\|y\|^p)$ and $L=2^{p-1}$ 
such that $\theta \geq0$ and $p\in[0,1)$, and let $\mathcal A$ be a normed algebra, then we get the following result.
\begin{corollary}\label{cor1}
Let $d:\mathcal A\times \mathcal A\to\mathcal X_{\rho}$ be a maaping with
$d(x,0)=d(0,x)=0$ and satisfy the inequality of the form
\begin{align*}
\rho\Big(d(\lambda (x+y),z+w)&+d(\lambda (x+y),z-w)+d(\lambda (x-y),z+w)\nonumber\\&+d(\lambda (x-y),z-w)-4\lambda d(x,z)\Big)\nonumber\\
&\leq \rho\Big(4s\big[d(\frac{x+y}{2},z-w)+d(\frac{x-y}{2},z+w)-d( x,z)+d(y,w)\big]\Big)\nonumber\\&+\theta(\|x\|^p+\|y\|^p)(\|z\|^p+\|w\|^p)
\end{align*}
for all $\lambda \in \mathbb T,p\in[0,1)$ and $x,y,z,w \in \mathcal A$. Then there
exists a unique bi-additive mapping $D:\mathcal A\times \mathcal A\to \mathcal X_{\rho}$ which is $\mathbb C$-linear in the first variable, such that
\begin{equation*}
\rho(D(x,z)-d(x,z))\leq \frac{\theta}{2^{1-p}-1}\|x\|^p\|z\|^p
\end{equation*}
for all $x,z \in \mathcal A$. In addition, the mapping $d:\mathcal A\times \mathcal A\to\mathcal X_{\rho}$ satisfies $d(2x,z)=2d(x,z)$ and
\begin{equation*}
\rho \big(d(xy,z)-d(x,z)y-xd(y,z)\big)\leq \psi(x,y)\psi(z,w),
\end{equation*}
for all $x,y,z,w \in \mathcal A$, then $d$ is a biderivation.
\end{corollary}
\begin{theorem}\label{th.2}
Let $d:\mathcal A\times \mathcal A\to\mathcal X_{\rho}$ satisfy in 
$d(x,0)=d(0,x)=0$ and the inequality of the form
\begin{align}\label{eq2.1}
\rho\Big(d(\lambda (x+y),z+w)&+d(\lambda (x+y),z-w)+d(\lambda (x-y),z+w)\nonumber\\&+d(\lambda (x-y),z-w)-4\lambda d(x,z)\Big)\nonumber\\
&\leq \rho\Big(4s\big[d(\frac{x+y}{2},z-w)+d(\frac{x-y}{2},z+w)-d( x,z)+d(y,w)\big]\Big)\nonumber\\&+\psi(x,y)\psi(z,w)
\end{align}
for all $\lambda \in \mathbb T, s<1$ and $x,y,z,w \in \mathcal A$, where $\psi:\mathcal
A^2\to [0, \infty)$
is a given mapping such that
\begin{equation*}
\psi(x,x)\leq \frac{L}{2} \psi(2x,2x)
\end{equation*}
and
\begin{equation}\label{eq2.2}
\lim_{n\to \infty}2^n\psi(\frac{x}{2^n},\frac{x}{2^n})=0,
\end{equation}
for all $x,z \in \mathcal A$ and a constant $0<L<1$. Then there
exists a unique bi-additive mapping $D:\mathcal A\times \mathcal A\to \mathcal X_{\rho}$ which is $\mathbb C$-linear in the first variable, such that
\begin{equation}\label{eq2.3}
\rho(D(x,z)-d(x,z))\leq \frac{1}{2(1-L)}\psi(x,x)\psi(z,0)
\end{equation}
for all $x,z \in \mathcal A$.

In addition, the mapping $d:\mathcal A\times \mathcal A\to\mathcal X_{\rho}$ satisfies $d(2x,z)=2d(x,z)$ and
\begin{equation}\label{eq2.4}
\rho\big(d(xy,z)-d(x,z)y-xd(y,z)\big)\leq \psi(x,y)\psi(z,w),
\end{equation}
for all $x,y,z\in \mathcal A$, then $d$ is a biderivation.
\begin{proof}
Let $\mathcal E$ and $\widetilde{\rho}$ be the same which are defined in the proof of Theorem \ref{th.1}. We define the function $\mathcal{T}:\mathcal
E_{\widetilde{\rho}}\to\mathcal E_{\widetilde{\rho}}$ as follows
\begin{equation*}
\mathcal{T}\delta(x,z):=2\delta(\frac{x}{2},\frac{z}{2})
\end{equation*}

for all $\delta \in \mathcal E_{\widetilde{\rho}}$. Let $\delta,
\gamma \in \mathcal E_{\widetilde{\rho}}$ and let $c\in[0,\infty]$
be an arbitrary constant with $\widetilde{\rho}(\delta-\gamma)\leq
c$. We have
\begin{equation*}
\rho(\delta(x,z)-\gamma(x,z))\leq c \psi(x,x)\psi(z,0)
\end{equation*}
for all $x,z\in \mathcal A$. By the assumption and the last
inequality, we get
\begin{equation*}
\begin{split}
\rho\left(2\delta(\frac{x}{2},\frac{z}{2})-2\gamma(\frac{x}{2},\frac{z}{2})\right)\leq
\kappa \rho\left(\delta(\frac{x}{2},\frac{z}{2})-\gamma(\frac{x}{2},\frac{z}{2})\right)\leq \kappa c \psi(\frac{x}{2},\frac{x}{2})\psi(\frac{z}{2},0)\leq \frac{\kappa L}{2}c \psi(x,x)\psi(z,0)
\end{split}
\end{equation*}
for all $x,z\in \mathcal A$. Hence,
$\widetilde{\rho}(\mathcal{T}\delta-\mathcal{T}\gamma)\leq
\frac{\kappa L}{2}\widetilde{\rho}(\delta-\gamma)$, for all $\delta, \gamma \in
\mathcal E_{\widetilde{\rho}}$ so $\mathcal{T}$ is a
$\widetilde{\rho}$--strict contraction.

Letting $x=y=\frac{x}{2},\lambda=1$ and $w=0$ in (\ref{eq2.1}), and since $\rho(.x)$ is increasing we get
\begin{equation}\label{eq2.6}
\rho\left(d(x,z)-2d(\frac{x}{2},z)\right)\leq\rho\left(2d(x,z)-4d(\frac{x}{2},z)\right)\leq \psi(\frac{x}{2},\frac{x}{2})\psi(z,0)
\end{equation}
for all $x,z\in\mathcal A$. Replacing $x$ by $\frac{x}{2}$ in (\ref{eq2.6}) we
have
\begin{equation*}\label{eq2.7}
\rho\Big(d(\frac{x}{2},z)-2d(\frac{x}{4},z)\Big)\leq \psi(\frac{x}{4},\frac{x}{4})\psi(z,0)
\end{equation*}
for all $x,z\in\mathcal A$. Since $\rho$ is convex modular and
satisfies the $\Delta_2$--condition, for all $x,z\in\mathcal A$ we
have
\begin{equation*}
\begin{split}
\rho\left(2d(\frac{x}{4},z)-\frac{1}{2}d(x,z) \right)&\leq  \frac{1}{2}\rho\big(4d(\frac{x}{4},z)-2d(\frac{x}{2},z)\big)+\frac{1}{2}\rho\big(2d(\frac{x}{2},z)-d(x,z)\big)\\
&\leq \frac{\kappa}{2}\psi(\frac{x}{4},\frac{x}{4})\psi(z,0)+\frac{1}{2}\psi(\frac{x}{2},\frac{x}{2})\psi(z,0).
\end{split}
\end{equation*}
Consequently 
\begin{equation*}
\begin{split}
\rho\left(2^2d(\frac{x}{2^2},z)-d(x,z)\right)&\leq  \kappa\left[\rho\big(2d(\frac{x}{4},z)-\frac{1}{2}d(x,z)\big)\right]\\
&\leq \kappa\left[\frac{\kappa}{2}\psi(\frac{x}{4},\frac{x}{4})\psi(z,0)+\frac{1}{2}\psi(\frac{x}{2},\frac{x}{2})\psi(z,0)\right]\\
&\leq \frac{\kappa ^2}{2}\psi(\frac{x}{4},\frac{x}{4})\psi(z,0)+\frac{\kappa}{2}\psi(\frac{x}{2},\frac{x}{2})\psi(z,0).
\end{split}
\end{equation*}
for all $x,z\in\mathcal A$. By induction we obtain
\begin{align}\label{eq2.8}
\rho\left(2^nd(\frac{x}{2^n},z)-d(x,z)\right)&\leq
\sum_{i=2}^n\frac{\kappa ^n}{2^{n-i+1}} \psi(\frac{x}{2^i},\frac{x}{2^i})\psi(z,0)+\frac{\kappa ^{n-1}}{2^{n-1}}\psi(\frac{x}{2},\frac{x}{2})\psi(z,0)\nonumber\\
&\leq \sum_{i=2}^n\frac{2^n}{2^{n-i+1}} \psi(\frac{x}{2^i},\frac{x}{2^i})\psi(z,0)+\frac{2^{n-1}}{2^{n-1}}\psi(\frac{x}{2},\frac{x}{2})\psi(z,0)\nonumber\\
&\leq \sum_{i=1}^n2^{i-1}\frac{L^i}{2^i} \psi(x,x)\psi(z,0)\nonumber\\
&\leq \frac{1}{2(1-L)}\psi(x,x)\psi(z,0)
\end{align}
for all $x,z\in\mathcal A$. 

The rest of the proof is similar to the proof of Theorem \ref{th.1}.
\end{proof}
\end{theorem}
As the same of Corollary \ref{cor1}, if  in Theorem \ref{th.2} we set $\psi(x,y)=\sqrt{\theta}(\|x\|^p+\|y\|^p)$ and $L=2^{p-1}$ 
such that $\theta \geq0$ and $p>1$, then we get the following result.
\begin{corollary}\label{cor2}
Suppose $d:\mathcal A\times \mathcal A\to\mathcal X_{\rho}$ with
$d(x,0)=d(0,x)=0$ and satisfies in (\ref{eq2.1}). Then there
exists a unique bi-additive mapping $D:\mathcal A\times \mathcal A\to \mathcal X_{\rho}$ which is $\mathbb C$-linear in the first variable, such that
\begin{equation*}
\rho(D(x,z)-d(x,z))\leq \frac{\theta}{2^{1-p}-1}\|x\|^p\|z\|^p
\end{equation*}
for all $x,z \in \mathcal A$.

In addition, the mapping $d:\mathcal A\times \mathcal A\to\mathcal X_{\rho}$ satisfies $d(2x,z)=2d(x,z)$ and
\begin{equation*}
\rho \big(d(xy,z)-d(x,z)y-xd(y,z)\big)\leq \psi(x,y)\psi(z,w),
\end{equation*}
for all $x,y,z\in \mathcal A$, then $d$ is a biderivation.
\end{corollary}
Now, we investigate the bi-additive $s$-functional inequality (\ref{eq00}).
\begin{lemma}\label{lem3}
Suppose $f:\mathcal A\times \mathcal A\to\mathcal X_{\rho}$ with
$f(x,0)=f(0,x)=0$ and satisfies in the inequality of the form
\begin{align}\label{eq3.1}
&\rho\Big(4\left[f(\lambda\frac{x+y}{2},z-w)+f(\lambda\frac{x-y}{2},z+w)-\lambda f(x,z)-\lambda f(y,w)\right]\Big)\nonumber\\
&\leq\rho\Big(s\left[f(x+y,z+w)+f(x+y,z-w)+f(x-y,z+w)+f(x-y,z-w)-4f(x,z)\right]\Big)
\end{align}
for all $\lambda \in \mathbb T,s<1$ and $x,y,z,w \in \mathcal A$. Then $f$ is biadditive.
\begin{proof}
Letting $x=y,w=0$ and $\lambda=1$ in (\ref{eq3.1}), we obtain $2f(2x,z)=4f(x,z)$, for all $x,z\in\mathcal A$. Since $|s|<1$, so we have
\begin{align}\label{eq3.2}
&\rho\left(f(x+y,z+w)+f(x+y,z-w)+f(x-y,z+w)+f(x-y,z-w)-2f(x,z)\right)\nonumber\\
&\leq \rho\left(4s\Big[f(\frac{x+y}{2},z-w)+f(\frac{x-y}{2},z+w)-f(x,z)+f(y,w)\Big]\right)\nonumber\\&=\rho\left(s\Big[2f(x+y,z-w)+2f(x-y,z+w)-4f(x,z)+4f(y,w)\Big]\right)
\end{align}
for all $x,y,z,w\in\mathcal A$. Letting $w=0$ in (\ref{eq3.1}), since $|s|<1$ we get 
$$f(x+y,z)+f(x-y,z)=2f(x,z),$$ and so $$f(x_1,z)+f(y_1,z)=2f(\frac{x_1+y_1}{2},z),$$ for all $x_1=x+y, y_1=x-y$ and $z\in \mathcal A$. Since $|s|<1$ and $f(0,z)=0$, for all $z\in \mathcal A$, we deduce that $f$ is additive in the first variable. Similarly, one can show that $f$ is additive in the second variable. Hence
$f$ is bi-additive.
\end{proof}
\end{lemma}
\begin{theorem}\label{th.4}
Let $d:\mathcal A\times \mathcal A\to\mathcal X_{\rho}$ satisfy in 
$d(x,0)=d(0,x)=0$ and the inequality of the form
\begin{align}\label{eq4.1}
&\rho\Big(4\big[d(\lambda\frac{x+y}{2},z-w)+d(\lambda\frac{x-y}{2},z+w)-\lambda d( x,z)-\lambda d(y,w)\big]\Big)\nonumber\\
&\leq\rho\Big(s\left[d(x+y,z+w)+d(x+y,z-w)+d(x-y,z+w)+d(x-y,z-w)-4d(x,z)\right]\Big)\nonumber\\&+\psi(x,y)\psi(z,w)
\end{align}
for all $\lambda \in \mathbb T, s<1$ and $x,y,z,w \in \mathcal A$, where $\psi:\mathcal
X^2\to [0, \infty)$
is a given mapping such that
\begin{equation*}
\psi(x,x)\leq \frac{L}{2}\psi(2x,2x)
\end{equation*}
and
\begin{equation}
\lim_{n\to \infty}2^n\psi(\frac{x}{2^n},\frac{x}{2^n})=0,
\end{equation}
for all $x,y \in \mathcal A$ and a constant $0<L<1$. Then there
exists a unique bi-additive mapping $D:\mathcal A\times \mathcal A\to \mathcal X_{\rho}$ which is $\mathbb C$-linear in the first variable, such that
\begin{equation}
\rho(D(x,z)-d(x,z))\leq \frac{1}{2(1-L)}\psi(x,x)\psi(z,0)
\end{equation}
for all $x,z \in \mathcal A$.

In addition, the mapping $d:\mathcal A\times \mathcal A\to\mathcal X_{\rho}$ satisfies $d(2x,z)=2d(x,z)$ and
\begin{equation}
\rho \big(d(xy,z)-d(x,z)y-xd(y,z)\big)\leq \psi(x,y)\psi(z,w),
\end{equation}
for all $x,y,z\in \mathcal A$, then $d$ is a biderivation.
\begin{proof}
Consider the set
$$\mathcal E=\{\delta:\mathcal A\times \mathcal A\to \mathcal X_\rho , \quad \delta(0,0)=0\}$$
and introduce the function $\widetilde{\rho}$ on $\mathcal E$ as
follows,
\begin{equation*}
\widetilde{\rho}(\delta)=\inf\{c>0 : \rho(\delta(x,y))\leq c \psi(x,0)\psi(y,0)\}.
\end{equation*}
As the proof of Theorem \ref{th.1} it is proved that $\widetilde{\rho}$ is convex modular and $\mathcal E_{\widetilde{\rho}}$ is $\widetilde{\rho}$-complete. Also the mapping $\mathcal{T}:\mathcal
E_{\widetilde{\rho}}\to\mathcal E_{\widetilde{\rho}}$ defined as 
\begin{equation*}
\mathcal{T}\delta(x,z):=\frac{1}{2}\delta(2x,2z)\quad (\delta\in\mathcal E_{\widetilde{\rho}}).
\end{equation*}
As proved in Theorem \ref{th.1}, $\mathcal T$ is a $\widetilde{\rho}$-strict contraction. 

Letting $y=w=0$ and $\lambda=1$ in (\ref{eq4.1}), and since $\rho(.x)$ is increasing we get
\begin{equation}\label{eq4.2}
\rho\left(2d(\frac{x}{2},z)-d(x,z)\right)\leq\rho\left(8d(\frac{x}{2},z)-4d(x,z)\right)\leq \psi(x,0)\psi(z,0)
\end{equation}
for all $x,z\in\mathcal A$. Replacing $x$ by $\frac{x}{2}$ in (\ref{eq4.2}) we
have
\begin{equation*}\label{eq4.3}
\rho\left(2d(\frac{x}{4},z)-d(\frac{x}{2},z)\right)\leq \psi(\frac{x}{2},0)\psi(z,0)
\end{equation*}
for all $x,z\in\mathcal A$. Since $\rho$ is convex modular and
satisfies the $\Delta_2$--condition, by induction we obtain
\begin{align}\label{eq4.4}
\rho\left(2^nd(\frac{x}{2^n},z)-d(x,z)\right)&\leq
\sum_{i=2}^n\frac{\kappa ^n}{2^{n-i+1}} \psi(\frac{x}{2^{i-1}},0)\psi(z,0)+\frac{\kappa ^{n-1}}{2^{n-1}}\psi(x,0)\psi(z,0)\nonumber\\
&\leq \sum_{i=2}^n\frac{2^n}{2^{n-i+1}} \psi(\frac{x}{2^{i-1}},0)\psi(z,0)+\frac{2^{n-1}}{2^{n-1}}\psi(x,0)\psi(z,0)\nonumber\\
&\leq \sum_{i=1}^n2^{i-1}\frac{L^i}{2^{i-1}} \psi(x,x)\psi(z,0)\nonumber\\
&\leq \frac{1}{2(1-L)}\psi(x,x)\psi(z,0)
\end{align}
for all $x,z\in\mathcal A$. The rest of the proof is similar to the proof of Theorem \ref{th.2}.
\end{proof}
\end{theorem}
if  in Theorem \ref{th.4} we set $\psi(x,y)=\sqrt{\theta}(\|x\|^p+\|y\|^p)$ and $L=2^{p-1}$ 
such that $\theta \geq0$ and $p>1$, then we get the following result.
\begin{corollary}\label{cor3}
Suppose $d:\mathcal A\times \mathcal A\to\mathcal X_{\rho}$ satisfies in  
$d(x,0)=d(0,x)=0$ and (\ref{eq4.1}). Then there
exists a unique bi-additive mapping $D:\mathcal A\times \mathcal A\to \mathcal X_{\rho}$ which is $\mathbb C$-linear in the first variable, such that
\begin{equation*}
\rho(D(x,z)-d(x,z))\leq \frac{\theta}{2^{1-p}-1}\|x\|^p\|z\|^p
\end{equation*}
for all $x,z \in \mathcal A$.

In addition, the mapping $d:\mathcal A\times \mathcal A\to\mathcal X_{\rho}$ satisfies $d(2x,z)=2d(x,z)$ and
\begin{equation*}
\rho \big(d(xy,z)-d(x,z)y-xd(y,z)\big)\leq \psi(x,y)\psi(z,w),
\end{equation*}
for all $x,y,z\in \mathcal A$, then $d$ is a biderivation.
\end{corollary}

\end{document}